\documentclass[12pt]{amsart}

\usepackage{type1cm}        
\usepackage{graphicx}        
\usepackage{multicol}        
\usepackage[bottom]{footmisc}
\usepackage[dvipsnames]{xcolor}

\usepackage{newtxtext}       %

\usepackage{amsmath,amssymb}   
\newcommand{\QQ}{{\mathbb Q}}

\newcommand{\RR}{{\mathbb R}}

\newcommand{\pp}{{\mathbb P}}
\newcommand{\CC}{{\mathbb C}}

\newcommand{\calX}{{\mathcal X}}

\title[Problems from Cortona]{Some problems at the interface of approximation theory\\ and algebraic geometry}

\author{F.~Sottile}
\address{Frank Sottile, Department of Mathematics,
         Texas A\&M University, College Station, Texas 77843,  USA}
\email{sottile@math.tamu.edu}
\urladdr{http://www.math.tamu.edu/\~{}sottile}
\thanks{This research stay was partially supported by National Science Foundation through grant  DMS-2201005 and the Simons
  Foundation and by the Mathematisches Forschungsinstitut Oberwolfach.} 

\begin{document}

\begin{abstract}
We sketch an assortment of problems that were posed---and not yet solved---during problem sessions at the conference
  ``Approximation Theory and Numerical Analysis meet Algebra, Geometry, and Topology'', which was held at the Palazzone Cortona 5--9
September 2022.
\end{abstract}

\maketitle


\section*{Introduction and Background}

Open problems are the bread and butter of mathematicians, and the congenial atmosphere at the Palazzone led several participants to share
some of theirs, and make comments about problems others posed.
This is an elaboration of a selection of those problems and comments, as understood by the author.
We try to give references and to follow the conventions established in the overview contributions to this volume~\cite{MST23,MS23}.
Any errors or misstatements are the fault of the author.


\section{Spline Spaces} \label{Sec:SplinesSpaces}

As explained in~\cite[Sect.\ 2.4]{MST23}, there are many open problems and directions of further study concerning the dimensions of spaces
of splines, even for planar domains.
Several specific problems along these lines arose during the problem sessions.

Given a triangulation or polyhedral subdivision $\Delta$ of a simply connected domain $D\subset\RR^2$ and positive integers $r,d$, we let
$S^r_d(\Delta)$ be the vector space of functions on $D$ that have $r$ continuous derivatives, and whose restriction to each two-dimensional
cell of $\Delta$ is a polynomial of degree at most $d$.
As Hal Schenck pointed out at the beginning of the problem sessions, the biggest open problem in this area, and a very hard one at that, is
to determine the dimensions of every spline space $S^r_d(\Delta)$.
That is, determine the Hilbert function of the spline module $S^r_\bullet(\Delta):=\bigoplus_{d\geq 0} S^r_d(\Delta)$.
We discuss three related problems that are less of a Moonshot.

\subsection*{The $S^1_3(\Delta)$ problem: Michael Di Pasquale}

\begin{figure}[htb]
  \centering
  \includegraphics{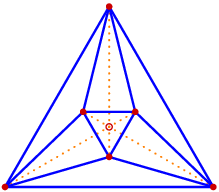}
   \qquad
  \includegraphics{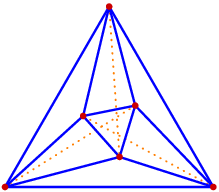}
  \caption{Morgan-Scott triangulation $\Delta$ on the left and a perturbation $\Delta'$  on the right.}
  \label{F:Morgan-Scott}
\end{figure}

This is a well-known outstanding problem in the theory of splines, which has been stated many times in the past.
For a planar triangulation $\Delta$, Schumaker~\cite{Schumaker79}  gave a lower bound for the dimension of a spline  space
$S^r_d(\Delta)$ that
depends only upon $f^\circ_1$, the number of interior edges in $\Delta$, $f^\circ_0$, the number of interior vertices, and the number of
distinct slopes at each vertex.
This was later shown to be the dimension of $S^r_d(\Delta)$ when $d\geq 3r+2$~\cite{AS87,Hong}.
When $r=1$ and $d\geq 4$, we have
 \begin{equation}\label{Eq:Schumaker}
    \dim S^1_d(\Delta)\ =\
      (1-f^\circ_0)\binom{d+2}{2} +
      f^\circ_1\binom{d}{2} + 3 f^\circ_0 + \sigma\,,
 \end{equation}
where $\sigma$ is the number of singular vertices, those interior vertices that have only two slopes among their incident edges.
For $d\geq 5$, this is due to Morgan and Scott~\cite{MoSc75}, and when $d=4$, this is due to Alfeld, Schumaker, and Piper~\cite{APS87}.
When $d=2$, equality does not hold in~\eqref{Eq:Schumaker} as  Morgan and Scott showed~\cite{MoSc77}:
Let $\Delta$ be the Schlegel diagram of a regular octahedron, displayed in Figure~\ref{F:Morgan-Scott}.
Then $S^1_2(\Delta)=7$, but if we perturb the configuration to a combinatorially equivalent one $\Delta'$ in which the dotted lines (images
of the three diagonals of the octahedron) do {\it not} meet, then $S^1_2(\Delta')=6$.
(This is also discussed in~\cite[Ex.\ 4]{MST23}.)

For $r=1$ and $d=3$, Schumaker's bound becomes $10-7 f^\circ_0 + 3 f^\circ_1+\sigma$, for any planar triangulation $\Delta$.
The challenge/problem is either to prove that this is the dimension of $S^1_3(\Delta)$ for any triangulation $\Delta$ of a simply connected 
planar domain $\Delta$, or to find a counterexample for which $S^1_3(\Delta)$ exceeds this number.
Any counterexample must also have unexpected splines in degree 2 \`a la Morgan-Scott.

\subsection*{Wang-Shi splits $WS_d$ of the triangle: Tom Lyche, Carla Manni, Hendrik Speleers}

This is a simple question about plane geometry which has consequences for the dimension of spline spaces of a particular
subdivision of a triangle that was studied in the proposers' paper~\cite[Sect.\ 2.2]{LMS22}.

Begin with a triangle and subdivide each edge into $d$ uniform segments.
These $3d$ segments have $3d$ endpoints along the perimeter of the triangle, which include the vertices.
Draw all $3d(d{-}1)$ lines connecting these points.
This gives a cross-cut partition of the triangle, called a \emph{Wang-Shi split}, $WS_d$, as this was originally proposed by Wang
and Shi~\cite{WS90}.

\begin{figure}[htb]
  \centering
 \includegraphics[height=138pt]{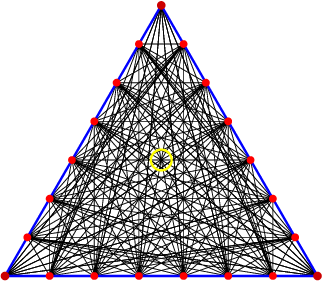}\quad
 \includegraphics[height=138pt]{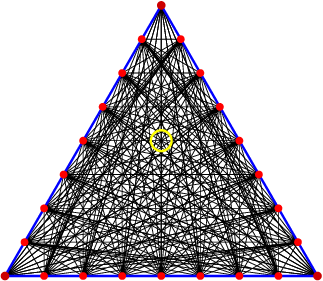}
 \caption{Wang-Shi splits $WS_7$ and $WS_8$ with vertices of maximal degrees circled.}
 \label{F:Wang-Shi}
\end{figure}

When $d=1$, this partition is just a triangle, and when $d=2$, this is the Powell-Sabin split~\cite{PS77}.
We display $WS_2$, $WS_3$, and $WS_4$.
\[
\includegraphics[height=80pt]{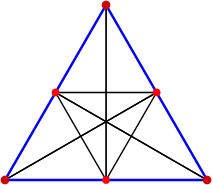}\quad
\includegraphics[height=80pt]{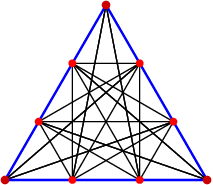}\quad
\includegraphics[height=80pt]{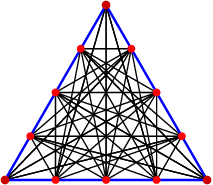}\quad
\]
Define  $\mu_d$ to be the maximum number of lines in $WS_d$ incident on a point in the interior of the triangle.
We have  $\mu_2=3$,  $\mu_3=3$, and $\mu_4=4$.
Table~\ref{tab:mu} displays values of $\mu_d$ for $d\leq 18$.
\begin{table}[htb] 
\caption{Values of $\mu_d$}
\label{tab:mu}      
%
%
\begin{tabular}{c c  c c c c c c r r r r r r r r r r }
$d$    &2&3&4&5&6&7&8&9&10&11&12&13&14&15&16&17&18\\\hline
$\mu_d$&3&3&4&5&6&7&{\color{blue}6}&8&8&9&9&10&10&10&12&{\color{blue}11}&12\rule{0pt}{9pt}\\
\end{tabular}
\end{table}
Figure~\ref{F:Wang-Shi} shows $WS_7$ and $WS_8$.
The circled vertices are those at which $\mu_7$ and $\mu_8$ achieve their maximum value.
The website\footnote{{\tt https://www.math.tamu.edu/\~{}sottile/research/stories/WangShi/index.html}} gives vector graphics pictures of
$WS_d$ for $d\leq 18$.
Note that $\mu_d$ is not a monotone function of $d$.
As the table shows, for $d$ equal to 8 and 17, $\mu_d < \mu_{d-1}$.

The problem is to determine $\mu_d$ for all $d$.
A subproblem is to show that $\mu_d\leq d{+}1$ for all $d$, which appears likely given the data.
Other problems include determining the total number of polygonal regions in $WS_d$ or the maximum number of sides of a polygon in $WS_d$
($WS_{15}$ has a decagon---a polygonal cell with ten sides).

The motivation  comes from the problem of determining the dimension of the spline spaces $S^r_k(WS_d)$.
The dimension of spline spaces for cross-cut partitions were determined by Chui and Wang~\cite{CW83}.
Their formula (given in Theorem 3.1 of~\cite{CW83}) is in terms of the number of cross-cuts (here, the number $3d(d-1)$ of lines), with a
correction term that is the sum over the internal vertices of a function of the number of lines incident at that vertex.

In~\cite[Sect.\ 2.2]{LMS22}, this formula is applied to the spline space $S^{d-1}_d(WS_d)$.
For this space, the Chui-Wang formula simplifies to become
\[
\dim S^{d-1}_d(WS_d)\ =\ \binom{d+2}{2} + 3d(d{-}1)\,,
\]
when  $\mu_d\leq d{+}1$.
Thus, the Wang-Shi splits are quite complicated polyhedral subdivisions which (conjecturally) support a relatively simple space of splines 
$S^{d-1}_d(WS_d)$. 

\subsection*{Overlap Splines: Oleg Davydov}

This problem arises in Davydov's contribution~\cite{OD} to this volume.
Fix positive integers $N,m,n$ and let $X \subset\RR^n$ be a set of $N$ points.
Suppose also that we have a system $\{X_i\subset X \mid i=1,\dotsc, m\}$ of subsets of $X$ whose union is $X$ (the sets form a cover of $X$).
Let us write $X$ for the set of points and the cover.

Let $P_d$ be the space of polynomials on $\RR^n$ of degree at most $d$.
An \emph{overlap spline of degree $d$ with respect to $X$} is a vector $s=\{p_i\in P_d\mid i=1,\dotsc,m\}$ of polynomials of degree at most
$d$ which satisfy the overlap condition
 \begin{equation}
 \label{Eq:overlap}
   p_i|_{X_i\cap X_j}\ =\  p_j|_{X_i\cap X_j} \qquad\mbox{for all } i,j\,.
 \end{equation}
(The condition is vacuous when $X_i\cap X_j=\emptyset$.)
Write $S_d(X)$ for the vector space of all such splines.
The problem is to determine the dimension of $S_d(X)$ in general.

Since  $S_d(X)$ is a linear subspace of $P_d^m$, we have the upper bound 
\[
  \dim(S_d(X))\ \leq\  m\cdot\dim(P_d)\ =\ m\binom{n+d}{n}\,.
\]
Rewriting the conditions~\eqref{Eq:overlap} in terms of $x\in X$  gives
\[
p_i(x)\ =\ p_j(x)\qquad \forall i,j\mbox{ with } x\in X_i\cap X_j\,,
\]
which are $m_x-1$ independent equations, where $m_x:=\#\{i\mid x\in X_i\}$.
As $\sum_{x\in X} m_x = \sum_{i=1}^m |X_i|$, we obtain the lower bound,
 \begin{equation}\label{Eq:lowerBound}
  \dim(S_d(X))\ \geq\  m\binom{n+d}{n} - \Bigl(\sum_{i=1}^m |X_i| -N\Bigr)\,,
 \end{equation}
which is an equality when the equations imposed by different $x\in X$ are independent.

There is one fairly straightforward case in which this is known, and it suggests some of the ingredients of an answer to this problem, or
some conditions that may be imposed on $X$ for a uniform solution. 
An \emph{interpolation set for $P_d$} is a finite set $Y\subset\RR^n$ of exactly $\dim(P_d)=\binom{n+d}{n}$ points which impose
independent conditions on polynomials in $P_d$.
Equivalently, the map from $P_d$ to $\RR^Y$ given by restricting polynomials in $P_d$ to functions on $Y$ is an isomorphism of vector spaces.
Proposition~1 in~\cite{OD} states that when  each $X_i$ is an interpolation set for $P_d$, then $\dim(S_d(X))=N$.
This follows from~\eqref{Eq:lowerBound} as $|X_i|=\binom{n+d}{n}$ when $X_i$ is an interpolation space for $P_d$.

An audience member mentioned that this has the flavor of \v{C}ech cohomology.
For each $1\leq i<j\leq m$, let $X_{i,j}:=X_i\cap X_j$ and write $\RR^{X_{i,j}}$ for the vector space of functions on $X_{i,j}$.
Then the overlap splines form the kernel of the map
 \begin{equation}\label{Eq:CechComplex}
   \bigoplus_{i=1}^m P_d\ \longrightarrow\ \bigoplus_{1\leq i<j\leq m} \RR^{X_{i,j}}\,,
 \end{equation}
in which a vector of polynomials $(p_i\mid i=1,\dotsc,m)$ is mapped to the vector
$( p_i|_{X_{i,j}}-p_j|_{X_{i,j}} \mid 1\leq i<j\leq m)$ of functions on the $X_{i,j}$.
The lower bound~\eqref{Eq:lowerBound} is obtained by analyzing the image of the map in~\eqref{Eq:CechComplex}.
This map~\eqref{Eq:CechComplex} has the flavor of the first coboundary operator in a \v{C}ech complex~\cite{BT}.

\section{Possible structures on splines}

Another class of problems concerned structures in spaces of splines inspired by other mathematical disciplines.

\subsection*{A Chinese Remainder Theorem for Splines?  Kiumars Kaveh}

This question is about extensions of splines.
Suppose that $\Delta$ is a triangulation (or subdivision) of a domain $D$ and that $D'\subset D$ is a simply connected domain
that is a union of some triangles or polyhedra in $\Delta$.
Write $\Delta'$ for the induced subdivision of $D'$, which we call a subtriangulation of $\Delta$.
Given a spline $g'$ on $\Delta'$, this question asks what are conditions on $g'$ so that there is a spline $g$ on $\Delta$  whose
restriction to $\Delta'$ is $g'$.
That is, under what conditions can a spline $g'$ on $\Delta'$  be extended to a spline $g$ on $\Delta$?

The motivation for the name of this question is that a spline is the solution to a system of congruences on vectors of polynomials.
It is not hard to find examples for which the map $S^r_d(\Delta)\to S^r_d(\Delta')$  is not surjective for any $d$.
Thus the pertinent question is to describe the image in general.
Di Pasquale's contribution to this volume~\cite{DiP} explains the lack of surjectivity and studies the image in many cases when the
extension of $\Delta'$ to $\Delta$ consists of adding a single simplex (in $\RR^2$, adding a single triangle).


More generally, we may have graph $G$ whose edges are labeled by ideals in a ring, or a GKM~\cite{GKM} graph  $G$, and $G'\subset G$ is an
induced subgraph.
We may again ask the same question about the surjectivity or iage of restriction.

An audience member asked if there is a version of Maier-Vietoris for splines.
That is, suppose that we have a triangulation or subdivision $\Delta$ of a domain that is the union of two subtriangulations $A,B$ whose
intersection $A\cap B$ is also a triangulation.
Then we have restriction maps
$S^r(\Delta)\to S^r(A)\oplus S^r(B) \to S^r(A\cap B)$, and a natural question would be under what conditions does this sequence have any
interesting homological properties.
This comes from the perspective that the space of splines on a subdivision is an analog of the global sections functor $H^0$ of that
subdivision.

\subsection*{Spline Schemes: Frank Sottile}

Let $\Delta$ be a polyhedral subdivision of a domain in $\RR^n$.
For any integer $r\geq 0$, the graded vector space $S^r_\bullet:=\bigoplus_{d\geq 0} S^r_d(\Delta)$ and the filtered vector space
$S^r(\Delta)=\bigcup_{d\geq 0} S^r_d(\Delta)$ are finitely generated modules over the rings $\RR[x_0,\dotsc,x_n]$ and $\RR[x_1,\dotsc,x_n]$,
respectively.
For $S^r_\bullet(\Delta)$, we use the homogenized versions of the spline spaces $S^r_d(\Delta)$, so that $S^r(\Delta)$ is a 
dehomogenization.

These modules of splines are are also finitely generated $\RR$-algebras.
To see that they are rings, let $I$ be an ideal and suppose that $(f_1,f_2)$ and $(g_1,g_2)$ are pairs whose differences lie in $I$.
That is, $f_1-f_2\in I$ and $g_1-g_2\in I$.
Then $(f_1g_1, f_2g_2)$ is another pair whose difference lies in $I$.
Indeed,
\[
  f_1g_1 - f_2g_2\ =\
  f_1g_1 -f_2g_1 + f_2g_1 - f_2g_2\ =\
  (f_1-f_2)g_1 + f_2(g_1-g_2)\ \in\ I\,.
\]
The spline ring is finitely generated over $\RR$ by the generators of the ring of constant splines (those that are global polynomial
functions, and thus take the same value at each $n$-dimensional cell in $\Delta$) and the module generators.

In algebraic geometry, the functor Proj associates a projective variety (scheme) to a graded algebra that is finitely generated in degree 1,
and Spec associates an affine variety (scheme) to a finitely generated algebra.
These may be embedded into  projective or affine space in such a way that the generators are the restrictions of coordinate functions.
For example, Spec$(\RR[x_1,\dotsc,x_N]/\langle f_1(x),\dotsc,f_m(x)\rangle$ is the set $\{x\mid f_i(x)=0,\ i=1,\dotsc,m\}$.
(Here $x\in\pp^{N-1}$ when the ring is graded, otherwise $x\in\RR^N$.)
Thus we may define the associated \emph{spline schemes}
\[
  \calX^r(\Delta)\ :=\ \mbox{Proj}(S^r_\bullet(\Delta))
  \qquad\mbox{and}\qquad
  X^r(\Delta)\ :=\ \mbox{Spec}(S^r(\Delta))\,.
\]
We now restrict to $S^r(\Delta)$ and $X^r(\Delta)$; while the same questions may be asked of $S^r_\bullet(\Delta)$ and $\calX^r(\Delta)$, we
expect they will have similar properties.
We will consider $X^r(\Delta)$ to be a subset of $\RR^N$, implicitly choosing generators of the spline module.

Let $\Delta_n$ be the set of $n$-dimensional cells in $\Delta$.
Then $S^{-1}(\Delta)$, the module/ring of discontinuous splines on $\Delta$, is
\[
 \RR[x_1,\dotsc,x_n]^{\Delta_n}\ :=\ \bigoplus_{\sigma\in\Delta_n} \RR[x_1,\dotsc,x_n]\,.
\]
Then $X^{-1}(\Delta)=\mbox{\rm Spec}(S^{-1}(\Delta))=\coprod_{\sigma\in\Delta_n}\RR^n$ is the disjoint union of affine spaces $\RR^n$, one
for each $n$-dimensional cell $\sigma$ of $\Delta$.

For each $-1\leq r<s$, we have inclusions of spline rings $S^s(\Delta)\hookrightarrow S^r(\Delta)$, which induce surjections of
spline schemes $X^r(\Delta)\twoheadrightarrow X^s(\Delta)$.
The map $X^{-1}(\Delta)\twoheadrightarrow X^r(\Delta)$ is a resolution of singularities of $X^r(\Delta)$,
and shows that $X^r(\Delta)$ has $|\Delta_n|$ irreducible components.
We use the adjective scheme for spline schemes, for in algebraic geometry, the term variety is often reserved for integral schemes, which
are reduced (no nilpotent elements of their coordinate rings) and irreducible.
Spline schemes are reduced (as spaces of splines are functions) but not irreducible.

The inclusion of $\RR[x_1,\dotsc,x_n]$ of the ring of constant splines induces a surjection
$\varphi\colon X^r(\Delta)\twoheadrightarrow\RR^n$, which
is a finite map as $S^r(\Delta)$ is a finite module over $\RR[x_1,\dotsc,x_n]$.
Also, while algebraic geometry often works best over $\CC$, the maps  $\coprod_{\sigma\in\Delta_n}\RR^n \to X^r(\Delta)\to \RR^n$ show that
there is no loss restricting to $\RR$.

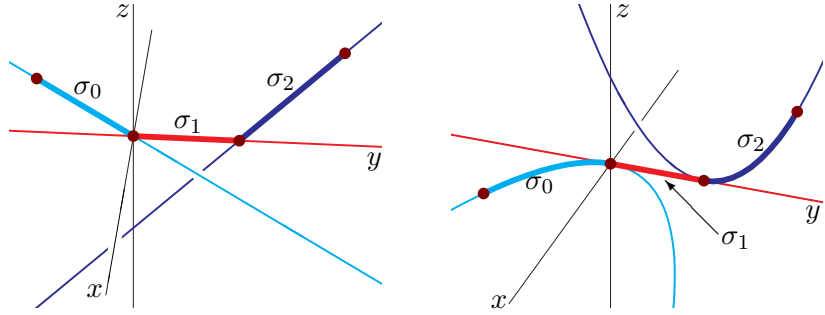
\begin{figure}[htb]
  \centering
  \begin{picture}(140,115)
      \put(0,0){\includegraphics[height=115pt]{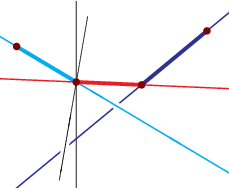}}
      \put(29,4){\small$x$}  \put(135,54){\small$y$}  \put(40,110){\small$z$}
      \put(24,81){\small$\sigma_0$} \put(62,68){\small$\sigma_1$} \put(97,84){\small$\sigma_2$}
   \end{picture}
  \qquad
    \begin{picture}(140,115)
      \put(0,0){\includegraphics[height=115pt]{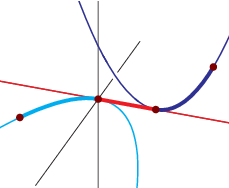}}
      \put(15,0){\small$x$}  \put(134,34){\small$y$}  \put(62,110){\small$z$}
      \put(28,45){\small$\sigma_0$} \put(102,24){\small$\sigma_1$} \put(108,61){\small$\sigma_2$}
                                    \put(101,28){\vector(-1,1){20}}
   \end{picture}
  \caption{Spline Schemes $X^0(\Delta)$ and $X^1(\Delta)$ for the subdivision $\Delta$ of~\eqref{Eq:Subdivision}.}
  \label{Fig:SplineSpaces}
\end{figure}

It is worthwhile to investigate these when $n=1$.
Consider the subdivision $\Delta$ of the interval $[-1,2]$ into three unit cells:
 \begin{equation}\label{Eq:Subdivision}
  \raisebox{-10pt}{\begin{picture}(130,22)(0,-10)

    \put(30,7){$\sigma_0$}    \put(60,7){$\sigma_1$}    \put(90,7){$\sigma_2$}
    \put(0,-1.5){\includegraphics{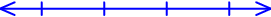}}
    \put(13,-10){$-1$}    \put(48,-10){$0$}    \put(78,-10){$1$}    \put(108,-10){$2$}

  \end{picture}}
 \end{equation}
Then we have
\[
S^0(\Delta)\ =\ \{(f_0(t),f_1(t),f_2(t))\in \RR[t]^3 \mid f_0(0)=f_1(0) \mbox{ and } f_1(1)=f_2(1)\}\,.
\]
The spline function represented by $(f_0,f_1,f_2)$ takes value $f_i(t)$ for $t\in\sigma_i$, and it is constant when $f_0=f_1=f_2$.
Observe that $S^0(\Delta)$ is generated as an $\RR[t]$-module by $(1,1,1)$, $x:=(t,0,0)$, and
$z:=(0,0,t-1)$.
If we let $y:=t(1,1,1)=(t,t,t)$, then $S^0(\Delta)$ is generated as a ring by $x,y,z$.
These satisfy the following relations $x(x-y)=xz=z(z-y{+}1))=0$, and thus
\[
S^0(\Delta)\ =\ \RR[x,y,z]/\langle x(x-y),xz,z(z-y{+}1) \rangle\,.
\]
We may check that
\[
\langle  x(x-y),xz,z(z-y{+}1)) \rangle\ =\
\langle x-y,z \rangle\ \cap\ 
\langle x,z \rangle\ \cap\ 
\langle x,z{-}y{+}1 \rangle\;,
\]
and thus $X^0(\Delta)$ is the union of three lines.
We display this on the left in Figure~\ref{Fig:SplineSpaces}.

Let us now consider $S^1(\Delta)$, which is
\[
\{(f_0(t),f_1(t),f_2(t))\in \RR[t]^3 \mid f_0-f_1\in\langle t^2\rangle \mbox{ and }
     f_1-f_2\in\langle (t-1)^2\rangle  \}\,.
\]
Observe that this is generated as a $\RR[t]$-module by $(1,1,1)$, $x:=(t^2,0,0)$ and $z:=(0,0,(t{-}1)^2)$.
If we let $y:=(t,t,t)$, then $S^1(\Delta)$ is generated as an $\RR$-algebra by $x,y,z$.
These satisfy the following relations $x^2-y^2x= xz= z^2- (y-1)^2z$, and thus
\[
S^1(\Delta)\ =\ \RR[x,y,z]/\langle x(x-y^2),xz,z(z-(y{-}1)^2) \rangle\,.
\]
We may check that
\[
\langle  x(x-y^2),xz,z(z-(y{-}1)^2)) \rangle\ =\
\langle x-y^2,z \rangle\ \cap\ 
\langle x,z \rangle\ \cap\ 
\langle x,z-(y{-}1)^2 \rangle\;,
\]
and thus $X^1(\Delta)$ is the union of the $y$-axis and two parabolas, $x=y^2$ in the $x,y$ plane and $z=(y-1)^2$ in the $y,z$ plane.
We display this on the right in Figure~\ref{Fig:SplineSpaces}.

Let $\sigma\in\Delta_n$ be a $n$-dimensional cell of $\Delta$.
Write $X^r_\sigma(\Delta)$ for the component of the spline scheme indexed by $\sigma$.
The composition $X^r_\sigma(\Delta)\hookrightarrow X^r(\Delta)\twoheadrightarrow \RR^n$ with the last map the finite
map $\varphi$ is an isomorphism.
Let $\widehat{\sigma}\subset X^r_\sigma(\Delta)$ be the inverse image of $\sigma$ under this isomorphism.
The union of these lifts over all $\sigma\in\Delta_n$ is a lift of $\Delta_n$ to $X^r(\Delta)$.
This forms the positive part $X^r_+(\Delta)$ of the spline scheme $X^r(\Delta)$.
In Figure~\ref{Fig:SplineSpaces}, these positive parts are indicated by the thickenings of the components.

It is not hard to generalize these examples to arbitrary subdivisions of $\RR^1$.
Spline schemes for arbitrary subdivisions in $\RR^d$ will have a strong combinatorial flavor.
One advantage of this perspective is that spline functions become restrictions of polynomial functions in $\RR^N$ to the
positive part of the spline scheme.

This can be asked for splines on graphs, or splines with varies smoothness conditions, such as supersmoothness~\cite{AS03}, as long as one
obtains a ring. 


\subsection*{Splines and the pushforward to a point: David Anderson}

As noted in~\cite[Sects.\ 1.2, 1.4]{MST23}, for some  algebraic manifolds $X$, the (torus) equivariant cohomology ring  $H^*_T(X)$ is 
equal to the space of  $C^0$-splines on a related subdivision $\Delta$ or labeled graph $G$.
This is the case when $X$ is a toric variety~\cite{Payne} or a GKM-space~\cite{Tymoczko}.
Equivariant cohomology has a natural map to the ring of polynomials, called integration or pushforward to a point.
For the spaces whose equivariant cohomology rings are also spaces of splines, the Atiyah-Bott formula~\cite[Cor.\ 5.2.4]{AndersonFulton} for
this pushforward is nicely expressed in terms of splines.
The question posed is to what extent can the Atiyah-Bott formula be generalized to arbitrary spline spaces?

We present two examples of this pushforward for toric varieties.
For these, the domain $D$ of the spline space is all of $\RR^n$ (or a restriction to a neighborhood of the origin), and the triangulation
$\Delta$ is a complete rational fan.
We begin with a one-dimensional example.
Let $\Delta$ be the union of two intervals $\sigma_-=(-\infty,0]$ and $\sigma_+=[0,\infty)$,
which form the 1-dimensional complete fan:
 \begin{equation}\label{Eq:1D_Fan}
  \raisebox{-10pt}{\begin{picture}(100,22)(0,-10)

    \put(23,7){$\sigma_-$}    \put(66,7){$\sigma_+$}  
    \put(0,-1.5){\includegraphics{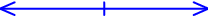}}
      \put(47.4,-10){$0$}  

  \end{picture}}
 \end{equation}
The corresponding toric variety is the projective line, $\pp^1$.
A $C^0$-spline (element of $H^*_T(\pp^1)$) is a pair $f=(f_-,f_+)$ of polynomials in $\RR[t]$  such that
$f_-(0)=f_+(0)$.
Here $f_-$ is the restriction of $f$ to $\sigma_-$, and  $f_+$ is the restriction of $f$ to $\sigma_+$.
Let $f=(-x^2+x, x^2-2x)$ be such a spline.
Then
\[
  \int f\ =\ \frac{-x^2+x}{-x} + \frac{x^2-2x}{x}\ =\
  (x-1) + (x-2)\  =\ 2x-3\,.
\]
(Notice that this is simply $(f_+-f_-)/x$.)

Let us now consider the projective plane $\pp^2$, which corresponds to the fan:
 \begin{equation}\label{Eq:2D_Fan}
  \raisebox{-42.5pt}{\begin{picture}(100,90)(-15,-5)

    \put(70,70){$\sigma_0$}  \put(57,-5){$\sigma_1$}   \put(-15,50){$\sigma_2$}  
 
     \put(0,0){\includegraphics{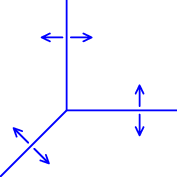}}
     \put(21,  0){$x-y$}    \put(-15, 26){$y-x$}
     \put(57, 14){$-y$}     \put( 65, 48){$y$}
     \put( 4, 65){$-x$}     \put( 47, 65){$x$}

  \end{picture}}
 \end{equation}
A $C^0$-spline is a triple $f=(f_0,f_1,f_2)$ of polynomials in $\RR[x,y]$  that satisfies $f_0-f_2\in\langle x\rangle$, 
$f_0-f_1\in\langle y\rangle$, and $f_1-f_2\in\langle x-y\rangle$.
The integral $\int f$ is defined to be the sum
\[
   \int f\ :=\ \frac{f_0}{xy}\ +\ \frac{f_1}{-y(x-y)}\ +\ \frac{f_2}{-x(y-x)}\,.
\]
The denominator of the summand corresponding to $\sigma_i$ is the product of the linear forms defining the boundary of $\sigma_i$, with signs
chosen so that they are positive on $\sigma_i$ (this is indicated in~\eqref{Eq:2D_Fan}).

For example, if $f=(xy,x^2-xy,-xy+y^2)$, then we have
\[
\int f\ =\ \frac{xy}{xy}\ +\ \frac{-xy+y^2}{-y(x-y)}\ +\ \frac{x^2-xy}{-x(y-x)}\ =\ 1+1+1\ =\ 3\,,
\]
and if $f=(x^2y+y^3, y^2x-x^2y, y^3-xy^2)$, then 
\[
\int f\ =\ \frac{x^2y+y^3}{xy}\ +\ \frac{y^2x-x^2y}{-y(x-y)}\ +\ \frac{y^3-xy^2}{-x(y-x)}\ =\ 2x\,.
\]

This integral makes sense for splines on any complete rational fan in $\RR^d$,
as such a spline is a class in equivariant cohomology.

We may use this to define the integral for splines on the decomposition $\Delta\subset\RR^1$ of the interval $[-1,2]$~\eqref{Eq:Subdivision}.
Indeed, let $f=(f_0,f_1,f_2)$ be a spline on $\Delta$.
Then
\[
\int f\ =\ \frac{f_0}{-t} + \frac{f_1}{t(1-t)} + \frac{f_2}{t-1}\ = \
\frac{f_1-f_0}{t} + \frac{f_2-f_1}{t-1}\,,
\]
as $\frac{f_1}{t(1-t)}=\frac{f_1}{t} + \frac{f_1}{1-t}$.
Note that not only is this well-defined, but if $f$ is homogeneous of degree $d$, then $\int f$ has degree $d{-}1$.

The question is if it is possible to generalize this pushforward to arbitrary splines, and if so, what are its properties?
For example, if $f$ is homogeneous of degree $d$, is $\int f$ a polynomial of degree $d{-}n$?
Some hurdles may be the choice of a normalization when the decomposition is not rational (defined over $\QQ$), and what to do if $\Delta$ is
not a triangulation.
We note that the pushforward to a point is an extremely important and natural map in equivariant cohomology.


\providecommand{\bysame}{\leavevmode\hbox to3em{\hrulefill}\thinspace}
\providecommand{\MR}{\relax\ifhmode\unskip\space\fi MR }
\providecommand{\MRhref}[2]{%
  \href{http://www.ams.org/mathscinet-getitem?mr=#1}{#2}
}
\providecommand{\href}[2]{#2}

\end{document}